\documentclass{amsart}
\usepackage{amsfonts}
\usepackage{amsmath}
\usepackage{amssymb}
\usepackage[utf8]{inputenc}

\newtheorem{theorem}{Theorem}
\newtheorem{corollary}{Corollary}
\newtheorem{lemma}{Lemma}
\newtheorem{proposition}{Proposition}
\newtheorem{notation}{Notation}
\newtheorem{definition}{Definition}
\newtheorem{remark}{Remark}
\newtheorem{example}{Example}

\begin{document}

\title[Graded identities of simple graded real division algebras]{Graded identities of simple real graded division algebras}

\author{Yuri Bahturin}
\address{Department of Mathematics and Statistics, Memorial University of Newfoundland, St. John's, NF, A1C5S7, Canada}
\email{bahturin@mun.ca}

\author{Diogo Diniz}
\address{Unidade Acadêmica de Matemática, Universidade Federal de Campina Grande, Campina Grande, PB, 58429-970, Brazil}
\email{diogo@mat.ufcg.edu.br}

\begin{abstract}
Let $A$ and $B$ be finite dimensional simple real algebras with division gradings by an abelian group $G$. In this paper we give necessary and sufficient conditions for the coincidence of the graded identities of $A$ and $B$. We also prove that every finite dimensional simple real algebra with a $G$-grading satisfies the same graded identities as a matrix algebra over an algebra $D$ with a division grading that is either a regular grading or a non-regular Pauli grading. Moreover we determine when the graded identities of two such algebras coincide. For graded simple algebras over an algebraically closed field it is known that two algebras satisfy the same graded identities if and only if they are isomorphic as graded algebras. 
\end{abstract}
\maketitle

\section{Introduction}

Two finite dimensional simple algebras over an algebraically closed field are  isomorphic if and only if they satisfy the same polynomial identities. This is a direct consequence of Amitsur-Levitzki theorem. Analogous results have been proved for Lie algebras by A. Kushkulei and Yu. Razmyslov \cite{KushkuleiRazmyslov}, for Jordan algebras by V. Drensky and M. Racine \cite{DrenskyRacine} and by E. Neher \cite{Neher} and for more general non-associative algebras by I. Shestakov and M. Zaicev \cite{ShestakovZaicev}. The analogous result for graded algebras was proved by P. Koshlukov and M. Zaicev \cite{PlamenZaicev} in the case of algebras graded by abelian groups and by E. Aljadeff and D. Haile \cite{AljadeffHaile} for arbitrary groups.

If the field is not algebraically closed, there exits non-isomorphic algebras that satisfy the same identities. For example, the real quaternion algebra and the algebra of $2\times 2$ matrices satisfy the same polynomial identities. The goal of this paper is to precisely determine, when two $G$-graded simple finite-dimensional algebras over the field of real numbers, satisfy the same $G$-graded identities, provided that $G$ is a finite abelian group.

For algebras over an algebraically closed field, the division gradings are regular (see Definition \ref{defreg}) and can be described in terms of a non-singular skew-symmetric bicharacter on the grading group (see \cite[Theorem 2.38]{ElduqueKochetov}). As it turns out, the graded polynomial identities for regular gradings are also determined by the corresponding bicharacter. Over the field of real numbers the situation is different. The division gradings for finite dimensional simple algebras over $\mathbb{R}$ have been classified in \cite{BahturinZaicev} and \cite{RodrigoEscudero}. In this classification, distinct regular gradings arise that have the same associated bicharacter and thus non-isomorphic algebras may have the same graded identities. Each division grading is obtained either by restriction of scalars in a complex algebra or as a tensor product of a non-regular and a regular component. In our main result, Theorem \ref{main}, we prove, roughly speaking, that the polynomial identities in each case are determined by the bicharacter of the regular component. In the last section we obtain similar results for finite dimensional simple real algebras. We prove in Theorem \ref{Thm2} that every finite dimensional graded simple real algebra with a $G$-grading satisfies the same graded identities as a matrix algebra over an algebra $D$ with a division grading that is either a regular grading or a Pauli grading. Then, in Theorem \ref{Thm3} we give necessary and sufficient conditions under which the graded identities of two such algebras coincide.

\section{Preliminaries}

We consider vector spaces, algebras, tensor products, etc. over the field $\mathbb{R}$ of the real numbers. Let $G$ be a group with identity element $e$ and $A$ an algebra. A vector space decomposition $A=\oplus_{g\in G} A_g$ is said to be a \textbf{$G$-grading} if $A_gA_h\subseteq A_{gh}$ for every $g,h \in G$. The \textbf{support} of the grading, denoted by $\mathrm{supp}\ A$, is the set $\mathrm{supp}\ A=\{g\in G| A_g\neq 0\}$. If $A=A_e$ we say that the grading is \textbf{trivial}. 

We say that two $G$-graded algebras $A$ and $B$ are \textbf{isomorphic} if there exists an isomorphism of algebras $\varphi:A\rightarrow B$ such that $\varphi(A_g)= B_g$. If $A$ is a $G$-graded algebra and $B$ is a $K$-graded algebra we say that $A$ and $B$ are \textbf{weakly isomorphic} if there exists an isomorphism of groups $\theta:G\rightarrow K$ and an isomorphism of algebras $\varphi:A \rightarrow B$ such that $\varphi(A_g)= B_{\theta(g)}$.

An element $a\in A$ is said to be \textbf{homogeneous} in the $G$-grading if there exists $g\in G$ such that $a\in A_g$. If $A$ has a unit and every non-zero homogeneous element is invertible we say that $A$ is a \textbf{graded division algebra}. 

\begin{example}
Let $G=(a)_2\times (b)_2 \cong \mathbb{Z}_2\times \mathbb{Z}_2$. Let $R=\mathbb{H}$ be the quaternion algebra with basis $\{1,i,j,k\}$ as a vector space and multiplication given by $i^{2}=j^{2}=-1$, $ij=-ji=k$. The decomposition \[R_e=\langle 1 \rangle, R_{a}=\langle i \rangle, R_{b}=\langle j \rangle, R_{c}=\langle k \rangle,\] where $c =ab$, is a $G$-grading on $\mathbb{H}$ that makes it a graded division algebra. We denote this graded division algebra by $\mathbb{H}^{(4)}$.
\end{example}

A division grading on $M_2(\mathbb{R})$ may be obtained by means of Sylvester matrices below: \[A=
\left(\begin{array}{cc}
	1 & 0\\
	0 & -1\\
\end{array}\right), B=
\left(\begin{array}{cc}
	0 & 1\\
	1 & 0\\
\end{array}\right), C=
\left(\begin{array}{cc}
	0 & 1\\
	-1 & 0\\
\end{array}\right).
\]

\begin{example}
Let $G=(a)_2\times (b)_2 \cong \mathbb{Z}_2\times \mathbb{Z}_2$. Let $S=M_2(\mathbb{R})$. The decomposition \[S_e=\langle I \rangle, R_{a}=\langle A \rangle, R_{b}=\langle B \rangle, R_{c}=\langle C \rangle,\] where $c=ab$ and $I$ is the identity matrix, is a $G$-grading on $M_2(\mathbb{R})$ that makes it a graded division algebra. We denote this graded division algebra by $M_2^{(4)}$.
\end{example}

The gradings in the previous two examples are regular in the sense that two homogeneous elements commute up to a non-zero real scalar. This generalized commutativity may be encoded in a function $\beta:G\times G \rightarrow \mathbb{R}^{\times}$, which is a \textbf{skew-symmetric bicharacter}, that is $\beta(g,h) \beta(h,g)=1$, for all $g,h \in G$ and also $\beta(gh,k)=\beta(g,k)\beta(h,k)$, for all $g,h,k \in G$. Next we recall the definition of a regular grading (see \cite{RegevSeeman}, \cite{AljadeffOfir}):
\begin{definition}\label{defreg}
Let $A$ be an associative algebra graded by a finite abelian group $G$. We say that the $G$-grading on $A$ is \textbf{regular} if there is a skew-symmetric bicharacter $\beta_A:G\times G\rightarrow \mathbb{R}^{\times}$ such that:
\begin{enumerate}
\item[(i)] For every integer $n\geq 1$ and every $n$-tuple $(g_1,g_2,\dots,g_n)\in G^{n}$, there are elements $a_i\in A_{g_i}$, $i=1,\dots,n$, such that $a_1\cdots a_n\neq 0$.
\item[(ii)] For every $g,h \in G$ and for every $a\in A_g$, $b\in A_h$ we have the equality \\
$ab=\beta_A(g,h)ba$.
\end{enumerate}
\end{definition}

Note that if $A$ has a regular $G$-grading then $\mathrm{supp}\ A=G$ by condition $(i)$. For a division grading we may consider $A$ graded by the subgroup $K=\mathrm{supp}\ A$. If the $K$-grading is regular then we say that $A$ has a regular grading by $\mathrm{supp}\ A$.

Over an algebraically closed field of characteristic zero every division $G$-grading on a finite dimensional simple algebra $A$ is regular (see \cite{BahturinZaicev2}). Thus if $A=M_n(\mathbb{C})$ is an algebra over $\mathbb{C}$ with a division grading, it is regular as a complex algebra with an associated complex bicharacter $\beta:G\times G\rightarrow \mathbb{C}^{\times}$. By restricting the scalars to $\mathbb{R}$ we obtain a real algebra with a division $G$-grading called \textbf{Pauli gradings} (see \cite{BahturinZaicev}). Note that $A$ is regular as a real algebra if and only if $\beta(G\times G)\subset \mathbb{R}^{\times}$. It follows from the classification of such gradings in \cite{BahturinZaicev2} that this is true if and only if $G$ is isomorphic to $\mathbb{Z}_2^{2k}$, for some $k$. Otherwise we have a non-regular grading on $A$ that we will refer to as a \textbf{non-regular Pauli grading}.

Now we give the definition of a graded polynomial identity. Let $f(x_{g_11},\dots, x_{g_nn})$ be a polynomial in the free associative algebra $\mathbb{R}\langle x_{gi}|g\in G, i \in \mathbb{N} \rangle$ and $A$ a $G$-graded algebra. A $n$-tuple $(a_1,\dots, a_n)$ of elements of $A$ is called a \textbf{$G$-admissible substitution} for the polynomial $f$ if $a_i\in A_{g_i}$, $i=i,\dots,n$. Moreover if $f(a_1,\dots,a_n)=0$ for every $G$-admissible substitution we say that $f$ is \textbf{graded polynomial identity} for the graded algebra $A$. We denote by $Id_G(A)$ the set of graded polynomial identities for $A$. If the group is trivial, $G=\{e\}$, we omit $G$ in the notation and write $Id(A)$ instead of $Id_G(A)$. Distinct letters, $x$, $y$, $z$, etc, may be used to denote distinct variables in a polynomial and in this case the second index $i$ will be omitted. If a polynomial $f(x_{g_11},\dots, x_{g_nn})$ is linear combination of monomials in which each of the variables $x_{g_11},\dots, x_{g_nn}$ appears exactly once we say that it is a \textbf{multilinear polynomial}.

\begin{remark}
It is well known that the set of graded polynomial identities of an algebra over a field of characteristic zero is determined by its multilinear elements. Hence to verify that $Id_G(A)=Id_G(B)$ for $G$-graded algebras $A$ and $B$ it is sufficient to verify that $A$ and $B$ have the same multilinear identities.
\end{remark}
It is clear that conditions $(i)$ and $(ii)$ of the Definition \ref{defreg} for the $G$-graded algebra $A$ may be written in terms of polynomial identities (see \cite[Lemma 28]{AljadeffOfir}): 
\begin{enumerate}
\item[(i)$^{\prime}$] $Id_G(A)$ contains no monomials with distinct indeterminates.
\item[(ii)$^{\prime}$] The polynomials $x_gy_h-\beta(g,h)y_hx_g$ are graded identities for $A$.
\end{enumerate}

Hence we obtain the following:

\begin{lemma}\label{typeI}
Let $A$ be an algebra with a regular $G$-grading and corresponding bicharacter $\beta_A$ and let $B$ be a $G$-graded algebra. Then $Id_G(A)=Id_G(B)$ if and only if  the above grading on $B$ is regular, with bicharacter $\beta_B=\beta_A$.
\end{lemma}

Thus in the case of regular division gradings, the polynomial identities are determined by the associated bicharacters. The algebras $\mathbb{H}^{(4)}$ and $M_2^{(4)}$ have a regular grading with the same bicharacter and by the previous lemma we have $Id_G(\mathbb{H}^{(4)})=Id_G(M_2^{(4)})$. However these algebras are not (graded) isomorphic. Next we give examples of graded division algebras whose gradings are not regular.

\begin{example}
Let $G=(a)_2\cong \mathbb{Z}_2$. The algebras $R=M_2(\mathbb{R})$ and $S=\mathbb{H}$ have the following division $G$-gradings:\[R_e=\langle I,C \rangle, R_{a}=\langle A, B \rangle\] and \[S_e=\langle 1,i \rangle, S_{a}=\langle j,k \rangle.\] These gradings will be denoted by $M_2^{(2)}$ and $\mathbb{H}^{(2)}$, respectively.
\end{example}

\begin{example}
Let $G=(a)_4\cong \mathbb{Z}_4$. A division $G$-grading on $R=M_{2}(\mathbb{C})$ is given by \[R_e=\langle I,C \rangle, R_{a}=\langle \omega A, \omega B \rangle,R_{a^2}=\langle i I, i C \rangle, R_{a^3}=\langle i \omega A, i\omega B \rangle,\] where $\omega = \frac{1}{\sqrt{2}}(1+i)$ is the $8$-th root of $1$ such that $\omega^2=i$. This grading will be denoted by $M_2(\mathbb{C}, \mathbb{Z}_4)$.
\end{example}

Let $K, L$ be two subgroups such that $G$ is the inner direct product $KL$. If $R$ is a $K$-graded algebra and $S$ is an $L$-graded algebra, the tensor product $R\otimes S$ of the algebras $R$ and $S$
has a $G$-grading, called the \textbf{tensor product grading}, which is given by $(R\otimes S)_g=R_k\otimes S_l$, where $g=kl$. In what follows $\mathrm{supp}\ R$ (resp. $\mathrm{supp}\ S$) will denote the support of the $K$-grading on $R$ (resp. the $L$-grading on $S$).

The tensor product of graded division algebras may not be a graded division algebra, for example, $\mathbb{C}\otimes \mathbb{C}$ is not a division algebra. Next we give two more examples of division gradings that appear in \cite{BahturinZaicev}.

\begin{example}
The tensor products $M_2^{(4)}\otimes \mathbb{C}^{(2)}$ and $\mathbb{H}^{(4)}\otimes \mathbb{H}$ are division gradings. Thus one obtains a $\mathbb{Z}_2^3$-divion grading on $M_2(\mathbb{C})$ denoted as $M_2^{(8)}$ and a $\mathbb{Z}_2^2$-division grading on $M_4(\mathbb{R})$, via the non-graded isomorphisms $\mathbb{H}\otimes \mathbb{H}\cong M_2(\mathbb{R})\otimes M_2(\mathbb{R})\cong M_4(\mathbb{R})$, denoted as $M_4^{(4)}$
\end{example}

Our study of the graded identities of graded division algebras relies on the classification given in \cite{BahturinZaicev}, their main result is the following:

\begin{theorem}\cite[Theorem 3.1]{BahturinZaicev}\label{BZ}
Any division grading on a real simple algebra $M_n(D)$, $D$ a real
division algebra, is weakly isomorphic to one of the following types

\begin{enumerate}
\item[$D=\mathbb{R}$:]
\begin{enumerate}
\item[(i)] $(M_2^{(4)})^{\otimes k}$;
\item[(ii)] $M_2^{(2)}\otimes(M_2^{(4)})^{\otimes (k-1)}$, a coarsening of $(i)$;
\item[(iii)] $M_4^{(4)}\otimes(M_2^{(4)})^{\otimes (k-2)}$, a coarsening of $(i)$
\end{enumerate}
\end{enumerate}

\begin{enumerate}
\item[$D=\mathbb{H}$:]
\begin{enumerate}
\item[(iv)] $\mathbb{H}^{(4)}\otimes(M_2^{(4)})^{\otimes k}$;
\item[(v)] $\mathbb{H}^{(2)}\otimes(M_2^{(4)})^{\otimes k}$, a coarsening of $(iv)$;
\item[(vi)] $\mathbb{H}\otimes(M_2^{(4)})^{\otimes k}$, a coarsening of $(v)$;
\end{enumerate}
\end{enumerate}

\begin{enumerate}
\item[$D=\mathbb{C}$:]
\begin{enumerate}
\item[(vii)] $\mathbb{C}^{(2)}\otimes(M_2^{(4)})^{\otimes k}$;
\item[(viii)] $\mathbb{C}^{(2)}\otimes M_2^{(2)}\otimes(M_2^{(4)})^{\otimes (k-1)}$, a coarsening of $(vii)$;
\item[(ix)] $(M_2^{(4)})^{\otimes k}\otimes \mathbb{C}^{(2)}\otimes \mathbb{H}$, a coarsening of $(vii)$; 
\item[(x)] $M_2^{(8)}\otimes(M_2^{(4)})^{\otimes (k-1)}$; 
\item[(xi)] $M_2(\mathbb{C},\mathbb{Z}_4)\otimes(M_2^{(4)})^{\otimes (k-1)}$, a coarsening of $(x)$;
\item[(xii)] $M_2^{(8)}\otimes(M_2^{(4)})^{\otimes (k-1)}\otimes \mathbb{H}$, a coarsening of $(x)$;
\item[(xiii)] Pauli grading.
\end{enumerate}
\end{enumerate}
None of the gradings of different types or of the same type but with different values
of $k$ is weakly isomorphic to the other.
\end{theorem}

\newpage
Next we list four types of gradings. 

\begin{enumerate}
\item[(I)] $A$ has a regular grading by $\mathrm{supp}\ A$ with bicharacter $\beta_A$;
\item[(II)] There exists a decomposition $G=K_AL_A$ of $G$ as a direct product of the subgroups $K_A$ and $L_A$ such that $A$ is isomorphic to $S_A\otimes R_A$ with a tensor product grading. The algebra $R_A$ has a regular grading such that $\dim (R_A)_e=1$ with support an elementary abelian 2 group $K_A$ and with bicharacter $\beta_{R_A}$. The algebra $S_A$ has grading with support $L_A$ that is weakly isomorphic to $\mathbb{H}^{(2)}$, $M_2^{(2)}$ or $M_2(\mathbb{C},\mathbb{Z}_4)$;
\item[(III)] $A$ is isomorphic to the algebra $\mathbb{H}\otimes R_A$ with a tensor product grading where $R_A$ has a regular $G$-grading such that $\dim (R_A)_e=1$ with bicharacter $\beta_{R_A}$ and $\mathbb{H}$ is trivially graded. 
\item[(IV)] $A$ has a non-regular Pauli grading with complex bicharacter $\beta_{A}$.
\end{enumerate}

It is clear that if two algebras have a regular grading their tensor product also has a regular grading. The algebras in the previous theorem have the following types:

\hspace{0.5cm}

\begin{center}
		\begin{tabular}{|c|c|}
		 \hline
			Type I & $(i)$, $(iv)$, $(vii)$, $(x)$ and regular Pauli gradings\\ \hline
			Type II& $(ii)$, $(v)$, $(viii)$ and $(xi)$\\ \hline
			Type III& $(iii)$, $(vi)$, $(ix)$ and $(xii)$\\ \hline
			Type IV& Non-regular Pauli gradings\\ \hline						
		\end{tabular}
\end{center}

\hspace{0.5cm}

Note that if the graded algebra $A$ is weakly isomorphic to the graded algebra $B$ then they have a grading of the same type $I-IV$. Hence it follows from the previous theorem that if $A$ is a finite dimensional simple real algebra with a division grading exactly one of the assertions $I-IV$ holds. This justifies the following:

\begin{definition}
We shall say that $A$ has a grading of Type $I-IV$ according to which assertion $I-IV$ above holds. For algebras of Type I and Type IV we let $R_A=A$, $\beta_{R_A}=\beta_{A}$.
\end{definition}

Now we are ready to state our main result:

\begin{theorem}\label{main}
Let $A$ and $B$ be two finite-dimensional simple $G$-graded division algebras. We have $Id_G(A)=Id_G(B)$ if and only if $A$ and $B$ have a grading of the same Type $(I)-(IV)$ above, $\mathrm{supp}\ A=\mathrm{supp}\ B$, $\mathrm{supp}\ R_A=\mathrm{supp}\ R_B$ and $\beta_A=\beta_B$ or $\beta_A=\overline{\beta_B}$.
\end{theorem}

\begin{remark}
Two graded algebras may satisfy the same graded identities and not be weakly isomorphic. The algebras $M_2(\mathbb{C})$ with a $\mathbb{Z}_2\times \mathbb{Z}_2$-Pauli grading and the  $\mathbb{Z}_2\times \mathbb{Z}_2$-graded algebra $M_2^{(4)}$ have regular gradings with the same bicharacter and thus satisfy the same graded identities. Note that these two algebras are not weakly isomorphic. 
\end{remark}

Theorem \ref{main} will be proved in Section \ref{proofmain}. We continue this section with some results that will be used in the proof.

Let $\theta:G\rightarrow K$ be a homomorphism of groups and $A$ a $G$-graded algebra. Denote by $\Gamma$ the $G$-grading $A=\oplus_{g\in G}A_g$ on $A$. The corresponding graded identities will also be denoted by $Id_G(\Gamma)$. We may obtain an $K$-grading on $A$, called \textbf{quotient grading induced by $\theta$} and denoted by $^{\theta}\Gamma$, as $^{\theta} \Gamma:A=\oplus_{k\in K}A_k$, where the homogeneous components are $A_k=\oplus_{g\in \theta^{-1}(k)}A_g$. The following result gives a useful criteria for the coincidence of the identities.

\begin{proposition}\label{propquot}
Let $A$ and $B$ be two algebras with $G$-gradings $\Gamma_A$ and $\Gamma_B$ respectively and $\theta:G\rightarrow K$ a homomorphism of groups. If $Id_G(\Gamma_A)=Id_G(\Gamma_B)$  then $Id_K(^{\theta}\Gamma_A)=Id_K(^{\theta}\Gamma_B)$.
\end{proposition}
\textit{Proof.}

Let $f(x_{k_11},\dots, x_{k_nn})$ be a $K$-graded identity for $A$. To prove that the polynomial $f$ is a $K$-graded identity for $B$ let $b_1,\dots, b_n$ be homogeneous elements of $B$ such that $b_i\in B_{k_i}$, $i=1,\dots,n$. We claim that $f(b_1,\dots, b_n)=0$. To prove this decompose $b_i$ into its homogeneous components in the $G$-grading \[b_i=b_{i1}+b_{i2}+\dots+b_{ik_i},\] for some integer $k_i$ and denote by $g_{ij}$ the homogeneous degree of $b_{ij}$. 
Note that since $b_i$ has homogeneous degree $k_i$ in the $K$-grading induced by $\theta$ we have $\theta(g_{ij})=k_i$. Let $S=\{(i,j)|1\leq i\leq n, 1\leq j \leq k_i\}$, $m=|S|$ and $\nu:S\rightarrow \{1,\dots, m\}$ a bijection. Construct the polynomial \[f^{\prime}=f\left(\sum_{j=1}^{k_1} x_{g_{1j}\nu((1,j))},\dots, \sum_{j=1}^{k_n} x_{g_{nj}\nu((n,j))}\right),\] obtained from $f$ by replacing $x_{k_i1}$ by the sum $\sum_{j=1}^{k_i} x_{g_{ij}\nu((i,j))}$. Since $\theta(g_{ij})=k_i$ and $f$ is a $K$ graded identity for $A$ we conclude that $f^{\prime}$ is a $G$-graded identity for $A$. Hence from the equality $Id_G(A)=Id_G(B)$ it follows that $f^{\prime}$ is a $G$-graded identity for $B$. Thus the result of the substitution $x_{g_{ij}\nu((i,j))}\mapsto b_{ij}$ is zero. Note that by the construction of $f^{\prime}$ the result of this substitution is $f(b_1,\dots, b_n)$. Hence $f(b_1,\dots, b_n)=0$ and $f$ is a $K$-graded identity for $B$. We proved the inclusion $Id_K(A)\subseteq Id_K(B)$. The reverse inclusion is proved analogously.
\hfill $\Box$

For finite dimensional graded simple algebras over an algebraically closed field two algebras satisfy the same graded identities if and only if they are isomorphic, see \cite{PlamenZaicev}, \cite{AljadeffHaile}. The next examples show that this is not true for algebras over the field of the real numbers.

\begin{example}\label{example1}
Let $G=(a)_2\times (b)_2\cong\mathbb{Z}_2\times \mathbb{Z}_2$, $K=(b)_2$ and $\beta$ be the bicharacter given by $\beta(g,h)=-1$ whenever $g\neq e$ and $h\neq e$. The algebras $M_2^{(4)}$ and $\mathbb{H}^{(4)}$ have a regular grading with skew-symmetric bicharacter $\beta$. Hence by Lemma \ref{typeI} we have $Id_{G}(M_2^{(4)})=Id_{G}(\mathbb{H}^{(4)})$. Moreover Proposition \ref{propquot} implies that $Id(M_2(\mathbb{R}))=Id(\mathbb{H})$.
\end{example}

\begin{example}\label{example2}
Let $G=(a)_2\times (b)_2$, $K=(b)_2$. Denote by $A$ the algebra $M_2(\mathbb{R})$ with $K$-grading such that the homogeneous component of degree $e$ consists of the diagonal matrices and the other homogeneous component consists of the off-diagonal matrices. This will be called the non-trivial elementary $G$-grading on $A$. Denote $\theta:G\rightarrow K$ be the homomorphism such that $\theta(a)=e$, $\theta(b)=b$. We may obtain $M_2^{(2)}$ and $A$ as quotient gradings (induced by $\theta$) of the algebras $R$ and $S$ respectively, where: \[R_{e}=\langle I \rangle, R_{a}=\langle A \rangle, R_{b}=\langle C \rangle, R_{c}=\langle B \rangle\] and \[S_{e}=\langle I \rangle, S_{a}=\langle B \rangle, S_{b}=\langle A \rangle, S_{c}=\langle C \rangle,\] where $c=ab$. Since $R$ and $S$ are regular gradings with the same skew-symmetric bicharacter $\beta$ as the algebra $M_2^{(4)}$ it follows from Lemma \ref{typeI} and Proposition \ref{propquot} that $Id_{K}(A)=Id_{K}(M_2^{(2)})$. Analogously it is proved that $Id_{K}(M_2^{(2)})=Id_{K}(\mathbb{H}^{(2)})$.
\end{example}

Our previous examples motivate the following definition:

\begin{definition}\label{def.}
Let $S,R$ be two algebras and identify $R$ (resp. $S$) with the subalgebra $1\otimes R$ (resp. $S\otimes 1$) of $S\otimes R$. A $G$-grading on $S\otimes R$ is said to be \textbf{standard} if the following two conditions hold:
\begin{enumerate}
\item[(i)] $S$ is a homogeneous subalgebra such that $Id_G(S)=Id_G(S^{\prime})$ for a matrix algebra $S^{\prime}$ with an elementary $G$-grading;
\item[(ii)] $R$ is a homogeneous graded division subalgebra such that $\dim R_e=1$ and $\mathrm{supp}\ R \cap \mathrm{supp}\ S=\{e\}$.
\end{enumerate}
\end{definition}

An algebra $A$ with a Type II grading has a standard grading if $S_A$ is weakly equivalent to $M_2^{(2)}$ or $\mathbb{H}^{(2)}$. In the next proposition we prove that for an algebra $A$, with a standard grading, the support of the regular component and its skew-symmetric bicharacter are determined by the identities of $A$. The proof is an adaptation of the arguments in \cite[Lemma 4]{PlamenZaicev}.

\begin{proposition}\label{prop2}
Let $A$ be a $G$-graded algebra isomorphic to the tensor product $S\otimes R$ of algebras $S,R$ with a standard grading such that $S_e$ is a commutative subalgebra. For every $g\in \mathrm{supp}\ A$ the polynomial $[x_e,y_g]$ is a graded identity for $A$ if and only if $g$ lies in $\mathrm{supp}\ R$. Moreover  if $\beta:K\times K \rightarrow \mathbb{R}^{\times}$ is the skew-symmetric bicharacter of $R$ then given $k_1,k_2\in K$ the polynomial $x_{k_1}y_{k_2}-\lambda y_{k_2}x_{k_1}$ is a graded identity for $A$ if and only if $\lambda = \beta(k_1,k_2)$.
\end{proposition}
\textit{Proof.}
For a matrix algebra $S^{\prime}$ with an elementary $G$-grading such that the neutral component is commutative, it is clear that $[x_e,y_{g^{\prime}}]$ is a graded identity if and only if $g^{\prime}=e$. Hence from condition $(i)$ of the previous definition we conclude that the same is true for $S$ as a $G$-graded algebra. Now let $K=\mathrm{supp}\ R$ and note that $\mathrm{supp}\ A=(\mathrm{supp}\ S)K$ and hence an element $g\in \mathrm{supp}\ A$ may be written as $g=g^{\prime}k$ where $g^{\prime}\in \mathrm{supp}\ S$ and $k\in K$. Fix an element $r\in R_k$, then it is clear that every element of $A_g$ may be written as a product $sr$ for some $s\in S_{g^{\prime}}$. Since $r$ is an invertible element that commutes with every element of $S$ the polynomial $[x_e,y_g]$ is a graded identity for $A$ if and only if $[x_e,y_{g^{\prime}}]$ is a graded identity for $S$. The last condition holds if and only if $g^{\prime}=e$ and the first assertion follows. 

Now we prove the second assertion. Let $r_1,r_2$ be non-zero elements of $R$ of degrees $k_1, k_2$, respectively. An element in $A_{k_i}$ has the form $a_ir_i$ where $a_i\in A_e$. Moreover $$(a_1r_1)(a_2r_2)-\lambda(a_2r_2)(a_1r_1)=(a_1a_2-\lambda \beta(k_2,k_1)a_2a_1)r_1r_2$$ and therefore $x_{k_1}y_{k_2}-\lambda y_{k_2}x_{k_1}$ is a graded identity for $A$ if and only if the polynomial $x_ey_e-\lambda \beta(k_2,k_1)y_ex_e$ is a graded identity for $S_e$. Since $x_ey_e-y_ex_e$ is also a graded identity for $S_e$ this is equivalent to the equality $\lambda \beta(k_2,k_1)=1$ and the result follows.
\hfill $\Box$

The remaining results in this section will be used to deal with Type II gradings that are not standard, i.e., for which $S_A$ is weakly equivalent to $M_2(\mathbb{C}, \mathbb{Z}_4)$.

\begin{proposition}\label{proptensor}
Let $A$ be a $G$-graded algebra and $C$ a commutative algebra with unit. If $A\otimes C$ is the tensor product algebra with the $G$-grading $(A\otimes C)_g=A_g\otimes C$ then $Id_G(A)=Id_G(A\otimes C)$. 
\end{proposition}
\textit{Proof.}
Let $f(x_{g_11},\dots, x_{g_nn})=\sum \lambda_{\sigma}x_{g_{\sigma(1)}\sigma(1)},\dots, x_{g_{\sigma(n)}\sigma(n)}$ be a multilinear polynomial. Note that $f$ is a graded identity for $A\otimes C$ if and only if it vanishes under admissible substitutions of the variables by elements of the form $a\otimes c$ where $a\in A$ is a homogeneous element and $c\in C$. Given $a_1,\dots, a_n \in A$ and $c_1,\dots, c_n \in C$ we have 
\begin{eqnarray*}
f(a_1\otimes c_1,\dots, a_n\otimes c_n)&=&\sum \lambda_{\sigma} a_{\sigma(1)}\otimes c_{\sigma(1)}\cdots a_{\sigma(n)}\otimes c_{\sigma(n)}\\&=&\sum \lambda_{\sigma} a_{\sigma(1)}\cdots a_{\sigma(n)}\otimes c_{\sigma(1)}\cdots c_{\sigma(n)}\\&=&\left(\sum \lambda_{\sigma} a_{\sigma(1)}\cdots a_{\sigma(n)}\right)\otimes c_{1}\cdots c_{n}\\&=&f(a_1,\dots, a_n)\otimes c_1\cdots c_n.
\end{eqnarray*}

Therefore a multilinear polynomial is an identity for $A$ if and only if it is an identity for $A\otimes C$ and this implies the equality $Id(A)=Id(A\otimes C)$.
\hfill $\Box$

\begin{example}\label{example}
Let $G= (a)_4 \cong \mathbb{Z}_4$, $G_2=\{g^2|g\in G\}$ and $\theta:G\rightarrow K$, where $K=G/G_2$, the quotient map. Consider the algebras $A=M_2(\mathbb{C})$ endowed  with the $G$-grading $\Gamma_A$ weakly isomorphic to $M_2(\mathbb{C},\mathbb{Z}_4)$ and $B=M_2(\mathbb{R})$ with the $K$-grading $\Gamma_B$ weakly isomorphic to $M_2^{(2)}$. It is clear that $^{\theta} \Gamma_A$ is isomorphic to $M_2^{(2)}\otimes \mathbb{C}$. Hence the equality $Id_{K}(^{\theta} \Gamma_A)=Id_{K}(\Gamma_B)$ follows from Proposition \ref{proptensor}.
\end{example}

\begin{notation}
Let $\theta:G\rightarrow K$ be a homomorphism of groups. Given a polynomial $f(x_{g_11},\dots, x_{g_nn}) \in \mathbb{R}\langle x_{gi}|g\in G, i \in \mathbb{N} \rangle$ we denote by $\theta(f)$ the polynomial $f(x_{\theta(g_1)1},\dots, x_{\theta(g_n)n}) \in \mathbb{R}\langle x_{ki}|k\in K, i \in \mathbb{N} \rangle$.
\end{notation}

\begin{proposition}\label{prop1}
Let $(g_1,\dots, g_n)$ be an $n$-tuple of elements of a group $G$ and let $f$ be a polynomial in $\mathbb{R}\langle x_{gi}|g\in G, i \in \mathbb{N} \rangle$ multilinear in the variables $x_{g_11},\dots, x_{g_nn}$. Moreover, let $A$ be an algebra with a $G$-grading $\Gamma_A$ and $\theta:G\rightarrow K$ a homomorphism of groups. A polynomial $\theta(f)$ is a $K$-graded identity for $^{\theta}\Gamma_A$ if and only if every polynomial in the set $S_{\theta,f}=\{f^{\prime}\in \mathbb{R}\langle x_{gi}|g\in G, i \in \mathbb{N} \rangle|\theta(f^{\prime})=\theta(f)\}$ is a $G$-graded identity for $A$.
\end{proposition}
\textit{Proof.}
Let $\textbf{B}$ be a basis for $A$ consisting of homogeneous elements. Note that $\theta(f)$ is a multilinear polynomial in $\mathbb{R}\langle x_{ki}|k\in K, i \in \mathbb{N} \rangle$, hence a graded identity if and only if it vanishes under every $K$-admissible substitution $(a_1,\dots, a_n)$ such that $a_i\in \textbf{B}$, $i=1,\dots, n$. If $b\in \textbf{B}$ is an element of $R_g$ for some $g\in G$ then $b\in (^{\theta}R)_{\theta(g)}$. Let $(b_1,\dots, b_n)$ be an $n$-tuple of elements of $\textbf{B}$ and let $g_i^{\prime}\in G$ such that $b_i\in B_{g_i^{\prime}}$, $i=1,\dots,n$. The tuple $(b_1,\dots, b_n)$ is a $K$-admissible substitution if and only if $\theta(g_i^{\prime})=\theta(g_i)$, $i=1,2,\dots, n$. Therefore $f$ is a $K$-identity if and only if it vanishes under every $\textbf{B}$-substitution $(b_1,\dots,b_n)$ such that $\theta(g_i^{\prime})=\theta (g_i)$, where $g_i^{\prime}$ is such that $b_i\in R_{g_i^{\prime}}$. This last assertion is equivalent to the claim that every element of $S_{\theta,f}$ is a $G$-graded identity for $A$. 
\hfill $\Box$

\begin{proposition}\label{pro2}
Let $G$ be an abelian group isomorphic to $\mathbb{Z}_4\times\mathbb{Z}_2^n$, $a\in G$ the generator of $G_2=\{g^2|g\in G\}$ and $\theta:G\rightarrow G/G_2$ the quotient map. Moreover, let $A=M_n(\mathbb{C})$ with a $G$-grading such that $iI\in A_{a}$ and let $f$, $f^{\prime}$ be polynomials in $\mathbb{R}\langle x_{gi}|g\in G, i \in \mathbb{N} \rangle$ which is multilinear in the variables $x_{g_11},\dots, x_{g_nn}$ and $x_{g_1^{\prime}1},\dots, x_{g_n^{\prime}n}$ respectively such that $\theta (f)=\theta (f^{\prime})$. In this case $f$ is a graded identity for $R$ if and only if $f^{\prime}$ is a graded identity for $R$.
\end{proposition}
\textit{Proof.}
Note that $\theta(f)=\theta(f^{\prime})$ if and only if $\theta(g_j)=\theta(g_j^{\prime})$, $j=1,2,\dots, n$. Now $\theta(g_j)=\theta(g_j^{\prime})$ if and only if $g_j^{\prime}=g_j$ or $g_j^{\prime}=a g_j$. 
We prove first that if $f$ is not a graded identity for $A$ then $f^{\prime}$ is not a graded identity for $A$. In this case there exists an admissible substitution $(r_1,\dots, r_n)$ such that $f(r_1,\dots, r_n)\neq 0$. Consider now the $n$-tuple $(r_1^{\prime},\dots, r_n^{\prime})$ where $r_j^{\prime}=r_j$ if $g_j^{\prime}=g_j$ and $r_j^{\prime}=ir_j$ if $g_j^{\prime}=a g_j$. It is clear that $(r_1^{\prime},\dots, r_n^{\prime})$ is an admissible substitution for $f^{\prime}$. Moreover $f^{\prime}(r_1^{\prime},\dots, r_n^{\prime})=i^{k}f(r_1,\dots, r_n)$, where $k$ is the number of $j$ such that $r_j^{\prime}=ir_j$. Therefore $f^{\prime}$ is not a graded identity. Analogously we can prove that if $f^{\prime}$ is not a graded identity then $f$ is not a graded identity. 
\hfill $\Box$

\begin{corollary}\label{z4}
Let $G$ be an abelian group isomorphic to $\mathbb{Z}_4\times\mathbb{Z}_2^n$, $a\in G$ the generator of $G_2=\{G^2|g\in G\}$, $\theta:G\rightarrow G/G_2$ the quotient map and $A=M_n(\mathbb{C})$ with a $G$-grading $\Gamma$ such that $iI$ is homogeneous of degree $a$. If $f$ is a polynomial in $\mathbb{R}\langle x_{gi}|g\in G, i \in \mathbb{N} \rangle$  multilinear in the variables $x_{g_11},\dots, x_{g_nn}$ then $\theta(f)$ is a graded identity for $^{\theta}\Gamma$ if and only if $f$ is a graded identity for $\Gamma$.
\end{corollary}
\textit{Proof.}
It follows from Proposition \ref{prop1} that $\theta(f)$ is a graded identity for $^{\theta}\Gamma$ if and only if every polynomial in the set $S_{\theta,f}=\{f^{\prime}\in \mathbb{R}\langle x_{gi}|g\in G, i \in \mathbb{N} \rangle|\theta(f^{\prime})=\theta(f)\}$ is a graded identity for $\Gamma$. By Proposition \ref{prop2} this last condition is equivalent to the condition that $f$ is a graded identity for $\Gamma$ and the result follows.
\hfill $\Box$

\section{Proof of the main result}\label{proofmain}

In this section we prove Theorem \ref{main}. The result will follow from Lemma \ref{typeI} together with the next three lemmas. 

\begin{remark}
Let $A$ be a $G$-graded algebra. Clearly $x_g\in Id_G(A)$ if and only if $g\notin \mathrm{supp}\ A$. Hence for $G$-graded algebras $A$ and $B$ the equality $Id_G(A)=Id_G(B)$ implies that $\mathrm{supp}\ A=\mathrm{supp}\ B$. Thus in the proof of Lemmas \ref{typeIII}, \ref{typeII} and \ref{typeIV} we omit the proof that $\mathrm{supp}\ A=\mathrm{supp}\ B$.
\end{remark}

\begin{lemma}\label{typeIII}
Let $A$ be an algebra with a $G$-grading of Type III with corresponding bicharacter $\beta_{R_A}$ and let $B$ be a $G$-graded algebra. Then $Id_G(A)=Id_G(B)$ if and only if $B$ has a $G$-grading of Type III, $\mathrm{supp}\ A=\mathrm{supp}\ B$ and $\beta_{R_A}=\beta_{R_B}$. 
\end{lemma}
\textit{Proof.}
The polynomial $x_ey_e-y_ex_e$ is a graded identity for every algebra with a grading of Types $I$, $II$, $IV$ and it is not a graded identity for any algebra with a grading of Type $III$. Since $Id_G(A)=Id_G(B)$ the algebra $B$ has a grading of Type $III$. The Hall polynomial $[[x_e,y_e]^2,z_e]$ is an identity for $A_e\cong \mathbb{H}$. We linearize on $y_e$ to obtain the identity \[[[x_e,y_e][x_e,w_e]+[x_e,w_e][x_e,y_e], z_e].\] Let $g,h \in G$ and fix $r\in (R_A)_g$ (resp. $r\in (R_A)_h$). Note that every element of $A_g$ (resp $A_h$) may be written as $ar$ (resp $bs$) where $a\in A_e$ (resp. $b \in A_e$). Hence \[([x_e,y_g][x_e,w_e]+[x_e,w_e][x_e,y_g]) z_h-\lambda z_h([x_e,y_g][x_e,w_e]+[x_e,w_e][x_e,y_g]) \] is a graded identity for $A$ if and only if $\lambda = \beta(g,h)$. This implies that if we have $Id_G(A)=Id_G(B)$ then $\beta_{R_A}=\beta_{R_B}$. Conversely if $B$ has a grading of Type III then the graded identities of $B$ (resp. $A$) are determined by $\mathrm{supp}\ B$ (resp. $\mathrm{supp}\ A$) and by the skew-symmetric bicharacter $\beta_{R_B}=\beta_{R_A}$ (see \cite[Theorem 3.4]{BahturinDrensky}) and hence $Id_G(B)=Id_G(A)$.
\hfill $\Box$

\begin{lemma}\label{typeII}
Let $A$ be an algebra with a $G$-grading of Type II with corresponding bicharacter $\beta_{R_A}$ and let $B$ be a $G$-graded algebra. Then $Id_G(A)=Id_G(B)$ if and only if $B$ has a $G$-grading of Type II, $\mathrm{supp}\ A=\mathrm{supp}\ B$, $\mathrm{supp}\ R_B=\mathrm{supp}\ R_A$ and $\beta_{R_A}=\beta_{R_B}$.
\end{lemma}
\textit{Proof.}
It follows from Lemmas \ref{typeI} and \ref{typeIII} that $B$ does not have a grading of Type I or of Type III respectively. Any algebra with a grading of Type IV satisfies the identities $[x_e,y_g]$ for every $g\in G$ and $A$ does not satisfy this set of identities since the algebra $S_A$ does not. Hence $B$ has a grading of Type II. In this case we have two possibilities: (1) $\mathrm{supp}\ R_A\cong (\mathbb{Z}_2)^n$ or (2) $\mathrm{supp}\ R_A\cong \mathbb{Z}_4\times (\mathbb{Z}_2)^n$. In the first case the grading is standard and it follows from Proposition \ref{prop2} that $\mathrm{supp}\ R_A=\mathrm{supp}\ R_B$ and $\beta_{R_A}=\beta_{R_B}$. In the second case let $G_2=\{g^2|g\in G\}$ and $\theta:G\rightarrow G/G_2$ the quotient map. We conclude from Example \ref{example} that the quotient grading is standard and therefore that $\mathrm{supp}\ R_A=\mathrm{supp}\ R_B$ and $\beta_{R_A}=\beta_{R_B}$ in this case.

Now we prove the converse. Assume that $B$ has a grading of Type II and that $\mathrm{supp}\ A =\mathrm{supp}\ B$, $\mathrm{supp}\ R_B=\mathrm{supp}\ R_A$ and $\beta_{R_A}=\beta_{R_B}$. We have two possibilities $\mathrm{supp}\ S_A\cong \mathbb{Z}_2$ or $\mathrm{supp}\ S_A\cong \mathbb{Z}_4$. Assume that $\mathrm{supp}\ S_A\cong \mathbb{Z}_2$. In this case $\mathrm{supp}\ A=\mathrm{supp}\ B$ is an elementary abelian $2$ group. Since $\mathrm{supp}\ B$ has no element of order $4$ we have $\mathrm{supp}\ S_B\cong \mathbb{Z}_2$. Moreover $S_A$ and $S_B$ are weakly isomorphic to one of the division algebras $\mathbb{H}^{(2)}$, $M_2^{(2)}$.

Let $K=\mathrm{supp}\ R_A=\mathrm{supp}\ R_B$, $L_A=\mathrm{supp}\ S_A=(a)_2$, $L_B=\mathrm{supp}\ S_B=(b)_2$. We denote by $S_A^{\prime}$ (resp. $S_B^{\prime}$) the algebra $M_2(\mathbb{R})$ with the non-trivial elementary $(a)_2$-grading (resp. $(b)_2$-grading). We denote by $A^{\prime}$ and $B^{\prime}$ the algebras $S_A^{\prime}\otimes R_A$ and $S_B^{\prime}\otimes R_B$, respectively, with the tensor product grading.  Since $Id_{L_A}(S_A)=Id_{L_A}(S_A^{\prime})$ and $Id_{L_B}(S_B)=Id_{L_B}(S_B^{\prime})$ (see Example \ref{example2}), it follows that $Id_{G}(A)=Id_{G}(A^{\prime})$ and $Id_{G}(B)=Id_{G}(B^{\prime})$ (see \cite[Theorem 3.4]{BahturinDrensky}).  The graded algebras $A^{\prime}$ and $B^{\prime}$ are isomorphic and hence $Id_G(A^{\prime})=Id_G(B^{\prime})$. Therefore we obtain the equality $Id_G(A)=Id_G(B)$. 
Now we consider the case where $\mathrm{supp}\ S_A\cong \mathbb{Z}_4$. In this case $\mathrm{supp}\ A=\mathrm{supp}\ B\cong \mathbb{Z}_4\times (\mathbb{Z}_2)^{2k}$ for some integer $k$. Hence we conclude that $\mathrm{supp}\ S_B\cong \mathbb{Z}_4$ and $S_A, S_B$ are weakly equivalent to the algebra $M_2(\mathbb{C},\mathbb{Z}_4)$. We denote $G^{\prime}=\mathrm{supp}\ A= \mathrm{supp}\ B$. Let $g_0$ be the generator of $G_2=\{g^2|g\in G^{\prime}\}$ and note that in this case $i1_A$ (resp. $i1_B$) is a homogeneous element of degree $g_0$. Let $\theta:G\rightarrow G/G_2$ be the canonical homomorphism. Let $\Gamma_A:A=\oplus_{g\in G}A_g$ and $\Gamma_B:B=\oplus_{g\in G}B_g$. Note that $^{\theta} \Gamma_A$ and 
$^{\theta}\Gamma_B$ have standard gradings (see Example \ref{example}). It follows from the previous case that $Id_{G/G_2}(^{\theta} \Gamma_A)=Id_{G/G_2}(^{\theta}\Gamma_B)$ and by Corollary \ref{z4} we conclude that $Id_G(A)=Id_G(B)$.
\hfill $\Box$

\begin{lemma}\label{typeIV}
Let $A$ be an algebra with a $G$-grading of Type IV and with corresponding complex bicharacter $\beta_{A}$ and let $B$ be a $G$-graded algebra. Then $Id_G(A)=Id_G(B)$ if and only if $B$ has a $G$-grading of Type IV with $\mathrm{supp}\ A= \mathrm{supp}\ B$ and with complex bicharacter $\beta_B$ such that either $\beta_A=\beta_B$ or $\beta_A=\overline{\beta_B}$.
\end{lemma}
\textit{Proof.}
We assume that $Id_G(A)=Id_G(B)$. It follows from Lemmas \ref{typeI}, \ref{typeIII} and \ref{typeII} that $B$ has a $G$-grading of Type IV. Let $\beta_B$ be its complex bicharacter. Given $g,h\in G$ either $\beta_A(g,h)\in \mathbb{R}$ or $\beta_A(g,h)$ is a complex root of a polynomial $t^2+qt+r$, for suitable real numbers $q,r$. Note that if $\beta_A(g,h)\in \mathbb{R}$ then \[x_{g}y_h-\beta_A(g,h)y_hx_g\] is a polynomial identity for $A$ and hence for $B$. In this case $\beta_B(g,h)=\beta_A(g,h)$. If $\beta_A(g,h)$ is a complex root of the polynomial $t^2+qt+r$, where $q,r \in \mathbb{R}$, then \[x_g^2y_h+qx_gy_hx_g+ry_hx_g^2\] is a graded identity for $A$ and hence for $B$. This implies that $\beta_B(g,h)$is also a root of the polynomial $t^2+qt+r$. Hence either $\beta_A(g,h)=\beta_B(g,h)$ or $\beta_A(g,h)=\overline{\beta_B(g,h)}$. In either case we conclude that given $g,h \in G$ either $\beta_A(g,h)=\beta_B(g,h)$ or $\beta_A(g,h)=\overline{\beta_B(g,h)}$. Now we prove that $\beta_A=\beta_B$ or $\beta_A=\overline{\beta_B}$.
The bicharacter $\beta_A$ induces a decomposition of the group $G$ as a direct product of cyclic subgroups: \[G=K_1^{\prime}\times K_1^{\prime \prime}\times \cdots \times K_r^{\prime}\times K_r^{\prime \prime}\] such that $K_i^{\prime}\times K_i^{\prime \prime}$ and $K_j^{\prime}\times K_j^{\prime \prime}$ are $\beta_A$-orthogonal for $i\neq j$, and $K_i^{\prime}$ and $K_i^{\prime \prime}$ are in duality by $\beta_A$. Note that since for every $g,h \in G$ we have $\beta_A(g,h)=\beta_B(g,h)$ or $\beta_A(g,h)=\overline{\beta_B(g,h)}$ the bicharacter $\beta_B$ induces the same decomposition above. Let $a_i$ be a generator of $K_i^{\prime}$ and $b_i$ be a generator of $K_i^{\prime \prime}$. If $\beta_A(a_i,b_i)=\beta_B(a_i,b_i)$ for $i=1,\dots, r$ then $\beta_A=\beta_B$. Otherwise we have $\beta_A(a_i,b_i)\neq\beta_B(a_i,b_i)$ for some $i$. We claim that $\beta_A=\overline{\beta_B}$. Assume on the contrary that $\beta_A(a_j,b_j)\neq \overline{\beta_B(a_j,b_j)}$ for some $j\neq i$. We proved that one of the two equalities
 $\beta_A(a_ia_j,b_ib_j)=\beta_B(a_ia_j,b_ib_j)$ or $\beta_A(a_ia_j,b_ib_j)=\overline{\beta_B(a_ia_j,b_ib_j)}$ holds. In the first case we obtain \[\beta_A(a_i,b_i)\beta_A(a_j,b_j)=\beta_B(a_i,b_i)\beta_B(a_j,b_j).\] Since $\beta_A(a_j,b_j)=\beta_B(a_j,b_j)$ this implies that $\beta_A(a_i,b_i)=\beta_B(a_i,b_i)$ which is a contradiction. In the second case we obtain \[\beta_A(a_i,b_i)\beta_A(a_j,b_j)=\overline{\beta_B(a_i,b_i)\beta_B(a_j,b_j)}.\] Since $\beta_A(a_i,b_i)=\overline{\beta_B(a_i,b_i)}$ we conclude that $\beta_A(a_j,b_j)=\overline{\beta_B(a_j,b_j)}$ which is also a contradiction. Clearly if $B$ has a $G$-grading of Type IV with complex bicharacter $\beta_B$ such that either $\beta_A=\beta_B$ or $\beta_A=\overline{\beta_B}$ then $Id_G(A)=Id_G(B)$ and the result is follows.
\hfill $\Box$

\section{Graded identities of finite-dimensional simple real algebras}

In this section $G$ denotes an abelian group. Let $(H, D, \textbf{g})$ be a triple where $H$ is a finite subgroup of $G$, $D$ a finite-dimensional simple real algebra with a division grading and support $H$ and  $\textbf{g}=(g_1,\dots,g_n)$ an $n$-tuple of elements of $G$. The algebra $R=M_n(D)$ has a $G$-grading $R=\oplus_{g\in G} R_g$ , determined by the triple $(H, D, \textbf{g})$, where \[R_g=\langle d_hE_{ij} |d_h \in D_h, g_i^{-1}hg_j=g \rangle.\]  We write $(H, D, \textbf{g}) \sim (H^{\prime}, D^{\prime}, \textbf{g}^{\prime})$ if $H=H^{\prime}$, $D\cong D^{\prime}$ and there exists $g\in G$ and a permutation $\pi \in S_n$ such that $g_iH=gg_{\pi(i)}^{\prime}H$, $i=1,\dots,n$.

\begin{remark}\label{isom}
A direct consequence of \cite{BahturinZaicev}[Theorem 7.3] is that any finite-dimensional simple real algebra with a $G$-grading is isomorphic to a graded algebra $M_n(D)$ determined by $(H,D,\textbf{g})$. Moreover two graded algebras $M_n(D)$ and $M_{n^{\prime}}(D^{\prime})$, determined by $(H,D,\textbf{g})$ and $(H^{\prime},D^{\prime},\textbf{g}^{\prime})$ respectively, are isomorphic if and only if $(H, D, \textbf{g}) \sim (H^{\prime}, D^{\prime}, \textbf{g}^{\prime})$.
\end{remark}

In this section we prove that any finite-dimensional simple real algebra with a grading by $G$ satisfies the same identities as an algebra $M_n(D)$, graded as above, where $D$ has either a division grading of Type I or Type IV. Moreover we determine when the graded identities of two such algebras coincide.

Next we consider the algebra $M_2(\mathbb{C},\mathbb{Z}_4)$. Let $G=(a)_4$, $H$ the subgroup generated by $a^2$ and $D$ the $H$-graded division algebra weakly isomorphic to $\mathbb{C}^{(2)}$. The $G$-grading on $R=M_2(\mathbb{C})$ determined by $(H,D,\textbf{g})$, where $\textbf{g}=(e,a)$, is given by
\begin{align}\label{grading}
 R_e=&
\left(\begin{array}{cc}
\mathbb{R} & 0\\
0 & \mathbb{R}
\end{array}\right), \hspace{0.3cm} R_a=
\left(\begin{array}{cc}
	0& \mathbb{R}\\
	i\mathbb{R}&0
\end{array}\right), 
\\ R_{a^2}=&
\left(\begin{array}{cc}
	i\mathbb{R} & 0\\
	0& i \mathbb{R}
\end{array}\right),
\hspace{0.3cm} R_{a^3}=
\left(\begin{array}{cc}
	0 & i\mathbb{R}\\
	\mathbb{R} & 0
\end{array}\right).\nonumber 
\end{align}

\begin{proposition}\label{pp1}
Let $G=(a)_4$. The algebra $M_2(\mathbb{C},\mathbb{Z}_4)$ satisfies the same $G$-graded identities as the algebra  $R=M_2(\mathbb{C})$ with the $G$-grading in (\ref{grading}).
\end{proposition}
\textit{Proof.}
Let $K=(b)_2$ and $\theta:G\rightarrow K$ the homomorphism such that $\theta(a)=b$. Denote by $\Gamma_1$ the division grading $M_2(\mathbb{C},\mathbb{Z}_4)$ and by $\Gamma_2$ the grading (\ref{grading}) on $M_2(\mathbb{C})$. Note that $^{\theta}\Gamma_1$ is isomorphic, as a graded algebra, to the tensor product $M_2^{(2)}\otimes \mathbb{C}$ and $^{\theta}\Gamma_2$ is isomorphic to the tensor product $R^{\prime}\otimes \mathbb{C}$ where $R^{\prime}$ denotes the algebra $M_2(\mathbb{R})$ with the non-trivial elementary $K$-grading. Since $Id_k(M_2^{(2)})=Id_k(R^{\prime})$ it follows that $Id_k(^{\theta}\Gamma_1)=Id_K(^{\theta}\Gamma_2)$. Therefore Corollary \ref{z4} implies that $Id_G(\Gamma_1)=Id_G(\Gamma_2)$.
\hfill $\Box$

\begin{theorem}\label{Thm2}
Every finite-dimensional simple real algebra with a grading by an abelian group $G$ satisfies the same $G$-graded identities as an algebra $M_n(D)$, with a grading determined by $(H, D, \textbf{g})$, where $D$ has a division grading of Type I or Type IV.
\end{theorem}
\textit{Proof.}
The result follows once we verify it holds for the algebras with  division gradings of Type II and Type III. For Type III gradings it follows from the equality $Id(\mathbb{H})=Id(M_2(\mathbb{R}))$. For Type II gradings the result follows from Example \ref{example2} and from the previous proposition.
\hfill $\Box$

If a graded division algebra $D$ has a regular grading or a non-regular Pauli grading the component $D_e$ consists of central elements, i.e., $D_e\subseteq Z(D)$. In the next lemma we prove that in this case $H$ is an invariant of the identities of $R=M_n(D)$. By Remark \ref{isom} the algebra $R$ is isomorphic to an algebra with a grading determined by $(H,D,\textbf{g})$ where $\textbf{g}=(g_1,\dots,g_n)$ is such that $g_i=g_j$ whenever $g_iH=g_jH$. In this case we have 
\begin{equation}\label{decomp}
R_e\cong M_{q_1}(D_e)\oplus\cdots\oplus M_{q_s}(D_e).
\end{equation}

\begin{remark}
In the proofs of Lemmas \ref{ll1} and \ref{ll2} we assume that $R$ is an algebra with a grading as above and that the decomposition (\ref{decomp}) holds.
\end{remark}

\begin{lemma}\label{ll1}
Let $R=M_n(D)$ be graded as in Theorem \ref{Thm2} and $f(x_{e1},\dots, x_{em})$ a multilinear central polynomial for $R_e$. For any $g\in G$ we have $g\in \mathrm{supp}\ D$ if and only if $[f,x_{g(m+1)}]\in Id_G(R)$ and $fx_{g(m+1)}\notin Id_G(R)$. 
\end{lemma}
\textit{Proof.}
Let $g \in \mathrm{supp}\ D$ and fix a nonzero element $d$ in $D_g$. An element $r$ in $R_g$ can be written as $r=dr^{\prime}$ for some $r^{\prime}\in R_e$. Given any $r_1,\dots, r_m \in R_e$ the element $r^{\prime \prime}=f(r_1,\dots, r_m)$ lies in $Z(R_e)$. Since $D_e\subseteq Z(R)$ we conclude that $[r^{\prime \prime}, r]=d[r^{\prime \prime},r^{\prime}]=0$. Hence $[f,x_{g(m+1)}]$ is a graded identity for $R$. The polynomial $f$ is a central polynomial for $R_e$ and hence not an identity. Therefore, we may assume that $r^{\prime \prime}\neq 0$. Moreover this element lies in $Z(R_e)$ and therefore $(r^{\prime \prime})^2\neq 0$. Let $d$ be a non-zero element in $D_g$. The element $s=dr^{\prime \prime}$ lies in $R_g$ and $r^{\prime \prime}s=d(r^{\prime \prime})^2\neq 0$. Hence $fx_{g(m+1)}$ is not a graded identity for $R$.

Now we prove the converse. Denote by $I_j$ the unit of the subalgebra of $R_e$ corresponding to $M_{q_j}$ in the decomposition (\ref{decomp}). We have $I=I_1+\cdots +I_s$, where $I$ is the unit of $R$. Since $fx_{g(m+1)}$ is not a graded identity for $R$, there exists an elementary matrix $E_{ij}$ in $R_{g^{\prime}}$, an element $d$ in $D_h$, with $g^{\prime}h = g$, and elements $r_1,\dots, r_m\in R_e$ such that $f(r_1,\dots, r_m)(dE_{ij})\neq 0$. We claim that $E_{ij}\in R_e$ and this proves the result. Indeed, in this case $g^{\prime}=e$ and $g=h\in \mathrm{supp} \ D$. It is clear that $E_{ij}=(I_1+\cdots+I_s)E_{ij}$ and hence $I_{l_1}E_{ij}\neq 0$ for some $l_1$. In this case, we have $I_{l_1}E_{ij}=E_{ij}$. Analogously, we have $E_{ij}I_{l_2}\neq 0$, for some $l_2$, and hence $E_{ij}I_{l_2}=E_{ij}$. The claim that $E_{ij}\in R_e$ follows from the equality $l_1=l_2$. Indeed, if we have $l_1=l_2$ then $E_{ij}=I_{l_1}E_{ij}I_{l_1}$ and this implies that $E_{ij}\in R_e$. Now we proceed to prove the equality $l_1=l_2$. The element $f(r_1,\dots,r_m)$ lies in the center of $R_e$ and therefore we may write $f(r_1,\dots,r_m)=d_1I_1+\cdots + d_sI_s$, where $d_1,\dots, d_s \in D_e$. Let $s_i$ denote the component of $r_i$ in the subalgebra of $R_e$ corresponding to $M_{q_{l_1}}$. It is clear that $f(s_1,\dots,s_m)=d_{l_1}I_{l_1}$. Note that $f(r_1,\dots, r_m)(dE_{ij})\neq 0$ and this implies that $d_{l_1}\neq 0$. Since $[f,x_{g(m+1)}]$ is a graded identity for $R$ we conclude that $[f(s_1,\dots,s_m),dE_{ij}]=0$ and therefore $E_{ij}I_{l_1}\neq 0$. Since $E_{ij}=E_{ij}I_{l_2}$ we have $(E_{ij}I_{l_2})I_{l_1}\neq 0$. This implies that $I_{l_2}I_{l_1}\neq 0$. Hence we conclude that $l_2=l_1$. 
\hfill $\Box$

Now we prove that the identities of the graded division algebra $D$ are determined by the graded identities of $R$.

\begin{lemma}\label{ll2}
Let $R=M_n(D)$ be an algebra with a grading as in Theorem \ref{Thm2} and let $f(x_{e1},\dots, x_{em},x_{e(m+1)})$ be a central polynomial for $R_e$ multilinear in $x_{e(m+1)}$. Denote by $f_{h_i}$ the polynomial obtained from $f$ by replacing the variable $x_{e(m+1)}$ by $x_{h_i(m+i)}$. The polynomial $g(x_{h_11},\dots, x_{h_kk})$ is an $H$-graded identity for $D$ if and only if $g(f_{h_1},\dots, f_{h_k})$ is a $G$-graded identity for $R$.
\end{lemma}

\textit{Proof.}
Let $g(x_{h_11},\dots, x_{h_kk})$ be an $H$-graded identity for $D$. Fix $r_1,\dots,r_m \in R_e$ and $r_{m+i}\in R_{h_i}$, $i=1,\dots,k$. We write $r_{m+i}=d_{h_i}r_{m+i}^{\prime}$, where $d_{h_i}\in D_{h_i}$ and $r_{m+i}^{\prime}\in R_e$. Since $f$ is a central polynomial for $R_e$ multilinear in $x_{e(m+1)}$ the element $\overline{f_{h_i}}=f_{h_i}(r_1,\dots,r_m,r_{(m+i)})$ can be written as $\overline{f_{h_i}}=d_{i1}I_1+\dots d_{is}I_s$, where $d_{ij}\in D_{h_i}$. The result of this substitution in $g(f_{h_1},\dots, f_{h_k})$ is \[g(\overline{f_{h_1}},\dots, \overline{f_{h_k})}=g(d_{11},\dots,d_{k1})I_1+\dots+g(d_{1s},\dots,d_{ks})I_s.\] Since $g$ is an $H$-graded polynomial identity for the algebra $D$ we conclude that $g(d_{1j},\dots,d_{kj})=0$ for $j=1,\dots,s$. Hence $g(\overline{f_{h_1}},\dots, \overline{f_{h_k})}$=0. This proves that $g(f_{h_1},\dots, f_{h_k})$ is a $G$-graded identity for $R$.

Now we prove the converse. The polynomial $f(x_{e1},\dots, x_{em},x_{e(m+1)})$ is a central polynomial for $R_e$ and hence it is not a graded identity. Let $r_1,\dots, r_m, r_{m+1}$ be elements of $R_e$ such that $f(r_1,\dots, r_m,r_{m+1})\neq 0$. It is clear that we may write $f(r_1,\dots, r_m,r_{m+1})=d_1^{\prime}I_1+\cdots+d_s^{\prime}I_s$, where $d_1^{\prime},\dots,d_s^{\prime} \in D_e$. There exists $i_0$ such that $d_{i_0}^{\prime}\neq 0$. Denote by $s_i$ the component of $r_i$ in the subalgebra of $R_e$ corresponding to $M_{q_{i_0}}$. Noting that $f(s_1,\dots, s_m,s_{m+1})=d_{i_0}^{\prime}I_{i_0}$, we replace $s_{m+1}$ by 
$t_{m+1}=(d_{i_0}^{\prime})^{-1}s_{m+1}$ to obtain $f(s_1,\dots, s_m,t_{m+1})=I_{i_0}$. Given the elements $d_1\in D_{h_1},\dots, d_k\in D_{h_k}$ we have $\overline{f_{h_i}}=f_{h_i}(s_1,\dots, s_m,d_it_{m+1})=d_iI_{i_0}$. Since $g(f_{h_1},\dots, f_{h_k})$ is a $G$-graded identity for $R$ we have \[0=g(\overline{f_{h_1}},\dots, \overline{f_{h_k}})=g(d_1,\dots, d_k)I_{i_0}.\] Therefore $g(d_1,\dots, d_k)=0$. This proves that $g$ is an $H$-graded identity for $D$.
\hfill $\Box$

\begin{proposition}\label{propos1}
Let $R=M_n(D)$ with $G$-grading $\Gamma$ determined by $(H,D,\textbf{g})$ and $\theta:G\rightarrow G/H$ the canonical homomorphism. The algebra $R$ with the quotient grading $^{\theta}\Gamma$ satisfies the same graded identities as the matrix algebra $M_{mn}(\mathbb{R})$ with the elementary grading induced by $(\theta(g_1),\dots,\theta(g_1),\dots,\theta(g_n))$, where $m$ is such that $Id(D)=Id(M_m(\mathbb{R}))$ and each entry $\theta(g_i)$ is repeated $m$ times. 
\end{proposition}
\textit{Proof.}
The grading $^{\theta}\Gamma:R=\oplus_{\overline{g}\in G/H}A_{\overline{g}}$ on $R$ is such that the homogeneous components are $R_{\overline{g}}=\langle dE_{ij}|d\in D,\theta(g_i)^{-1}\theta(g_j)=\overline{g}\rangle$. Since $D$ satisfies the same ordinary identities as $S=M_m(\mathbb{R})$ the algebra $R$ satisfy the same $G/H$-graded identities as $R^{\prime}=M_n(S)$ with the grading such that $R_{\overline{g}}^{\prime}=\langle sE_{ij}|s\in S,\theta(g_i)^{-1}\theta(g_j)=\overline{g}\rangle$. The canonical isomorphism of algebras $\varphi:M_n(S)\rightarrow M_{nm}(\mathbb{R})$ is a graded isomorphism if $M_{nm}(\mathbb{R})$ has an elementary grading determined by the $nm$-tuple $(\theta(g_1),\dots,\theta(g_1),\dots,\theta(g_n))$. Hence $R^{\prime}$ satisfy the same $G/H$-graded identities as $M_{nm}(\mathbb{R})$. 

\hfill $\Box$

\begin{lemma}\label{elem}
Let $\Gamma$ and $\Gamma^{\prime}$ be the elementary $G$-gradings on $M_n(\mathbb{R})$ induced by the $n$-tuples $\textbf{g}$ and $\textbf{g}^{\prime}$ respectively. The equality $Id_G(\Gamma)=Id_G(\Gamma^{\prime})$ occurs if and only if $g\in G$ and there exists a permutation $\pi \in S_n$ such that $g_iH=gg_{\pi(i)}H$, $i=1,\dots,n$.
\end{lemma}
\textit{Proof.}
It is clear that if the real $G$-graded algebras $R$ and 
$R^{\prime}$ satisfy the same graded identities then the complex algebras $R\otimes \mathbb{C}$ and $R^{\prime}\otimes \mathbb{C}$ satisfy the same graded identities as complex algebras. Hence the graded identities of the algebra $M_n(\mathbb{C})$ with elementary gradings induced by $\textbf{g}$ and $\textbf{g}^{\prime}$ coincide. Thus the result follows from \cite[Theorem 2]{PlamenZaicev}.
\hfill $\Box$

\begin{theorem}\label{Thm3}
Let $R=M_n(D)$ and $R^{\prime}=M_{n^{\prime}}(D^{\prime})$, with $G$-gradings determined by $(H, D, \textbf{g})$ and $(H^{\prime}, D^{\prime}, \textbf{g}^{\prime})$ respectively, where $D$ and $D^{\prime}$ have division gradings of Type I or Type IV. The algebras $R$ and $R^{\prime}$ satisfy the same $G$-graded identities if and only if $n=n^{\prime}$, $H=H^{\prime}$, $D$ and $D^{\prime}$ satisfy the same $H$-graded identities and there exists $g\in G$ and a permutation $\pi \in S_n$ such that $g_iH=gg_{\pi(i)}H$, $i=1,\dots,n$.
\end{theorem}
\textit{Proof.}
Assume first that $n=n^{\prime}$, $H=H^{\prime}$, $Id_H(D)=Id_H(D^{\prime})$ and that there exists $g\in G$ and a permutation $\pi \in S_n$ such that $g_iH=gg_{\pi(i)}H$, $i=1,\dots,n$. Denote $R^{\prime \prime}=M_n(D)$ with the $G$-grading determined by $(H,D,\textbf{g}^{\prime})$. Since $Id_H(D)=Id_H(D^{\prime})$ we have $Id_G(R^{\prime \prime})=Id_G(R^{\prime})$. Notice that $(H,D,\textbf{g}^{\prime}) \sim (H, D, \textbf{g})$ and therefore $R^{\prime \prime}$ is isomorphic as a graded algebra to $R$. Thus $Id_G(R)=Id_G(R^{\prime \prime})$. This proves that $Id_G(R)=Id_G(R^{\prime})$.

Now we prove the converse. It follows from Lemmas \ref{ll1} and \ref{ll2} that $H=H^{\prime}$ and that $D$ and $D^{\prime}$ satisfy the same graded identities. Note that $D$ (resp. $D^{\prime}$)  satisfies the same ordinary identities  as a matrix algebra $M_{m}(\mathbb{R})$ (resp. $M_{m^{\prime}}(\mathbb{R})$). It follows from Proposition \ref{propquot} that $Id(D)=Id(D^{\prime})$ and $Id(R)=Id(R^{\prime})$. The first equality implies that $m=m^{\prime}$ and the second that $nm=n^{\prime}m^{\prime}$. Thus we obtain $n=n^{\prime}$. It remains to prove the existence of the element $g$ and the permutation $\pi$. Let $\theta:G\rightarrow G/H$ denote the canonical homomorphism and $\Gamma, \Gamma^{\prime}$ the gradings on $R, R^{\prime}$ respectively. By Proposition \ref{propos1} we have that $Id_{G/K}(^{\theta}\Gamma)$ (resp.$Id_{G/K}(^{\theta}\Gamma^{\prime})$) is the set of identities of $M_{nm}(\mathbb{R})$ with the elementary grading induced by the $nm$-tuple $(\theta(g_1),\dots,\theta(g_1),\dots,\theta(g_n))$ (resp. $(\theta(g_1^{\prime}),\dots,\theta(g_1^{\prime}),\dots,\theta(g_n^{\prime}))$), where each $\theta(g_i)$ (resp. $\theta (g_i^{\prime})$) is repeated $m$ times. From Proposition \ref{propquot} we conclude that $Id_{G/K}(^{\theta}\Gamma^{\prime})=Id_{G/K}(^{\theta}\Gamma^{\prime})$. This and Lemma \ref{elem} imply that there exists $g\in G$ and $\pi^{\prime}\in S_{mn}$ such that $\theta(g_i)=\theta(gg_{\pi^{\prime}(i)}^{\prime})$, for $i=1,2,\dots, nm$. Hence the sets of entries in the $mn$-tuples $(\theta(g_1),\dots,\theta(g_1),\dots,\theta(g_n))$  and  $(\theta(gg_1^{\prime}),\dots,\theta(gg_1^{\prime}),\dots,\theta(gg_n^{\prime}))$ coincide and each entry appears the same number of times in both $mn$-tuples. Hence in the $n$-tuples $(\theta(g_1),\theta(g_2),\dots, \theta(g_n))$ and $(\theta(gg_1^{\prime}),\dots, \theta(gg_n^{\prime}))$ the sets of entries also coincide and each entry appears the same number of times in both tuples. Therefore there exists $\pi \in S_n$ such that $\theta(g_i)=\theta(gg_\pi(i)^{\prime})$, for $i=1,\dots, n$. 
\hfill $\Box$


\begin{thebibliography}{5}
\bibitem{AljadeffOfir} E. Aljadeff, D. Ofir, \textit{On regular G-gradings}, Trans. Amer. Math. Soc., \textbf{367} (2015), no. 6, 4207-4233.
\bibitem{AljadeffHaile} E. Aljadeff, D. Haile, \textit{Simple G-graded algebras and their polynomial identities}, Trans. Amer. Math. Soc., \textbf{366} (2014), no. 4, 1749--1771.
\bibitem{BahturinDrensky} Y. Bahturin, V. Drensky, \textit{Graded polynomial identities of matrices}, Linear Algebra and its Applications \textbf{357} (2002), 15--34.
\bibitem{BahturinZaicev} Y. Bahturin, M. Zaicev, \textit{Simple graded division algebras over the field of real numbers}, Linear Algebra and its Applications \textbf{490} (2016), 102--123.
\bibitem{BahturinZaicev2} Y. Bahturin, M. Zaicev, \textit{Graded algebras and graded identities}. Polynomial identities and combinatorial methods (Pantelleria, 2001) Lect. Notes in Pure and Appl. Math. \textbf{235} Dekker, New York, 2003, 101--139.
\bibitem{BahturinRegev} Y. Bahturin, A. Regev, \textit{Graded tensor products}, J. Pure Appl. Algebra \textbf{213} (2009), no. 9, 1643--1650.
\bibitem{DrenskyRacine} V. Drensky and M. Racine, \textit{Distinguishing simple Jordan algebras by means of polynomial identities}, Comm. Algebra \textbf{20} (1992), 309--327.
\bibitem{ElduqueKochetov} A. Elduque, M. Kochetov, \textit{Gradings on Simple Lie Algebras}, AMS Math. Surv. Monographs \textbf{189} (2013), 336p.
\bibitem{KushkuleiRazmyslov} A. Kushkulei, Yu. Razmyslov, \textit{Varieties generated by irreducible representations of Lie algebras}. (Russian), Vestnik Moskov. Univ. Ser. I Mat. Mekh., \textbf{5} (1983), 4--7.
\bibitem{Neher} E. Neher, \textit{Polynomial identities and nonidentities of split Jordan pairs}, Journal of Algebra \textbf{211} (1999), 206--224.
\bibitem{PlamenZaicev} P. Koshlukov, M. Zaicev, \textit{Identities and isomorphisms of graded simple algebras}, Linear Algebra and its Applications \textbf{432} (2010), 3141--3148.
\bibitem{RegevSeeman} A. Regev, T. Seeman, \textit{$\mathbb{Z}_2$-graded tensor products or p.i. algebras}, J. Algebra \textbf{291} (2005), no.1, 274--296.
\bibitem{RodrigoEscudero} A. Rodrigo-Escudero, \textit{Classification of division gradings on finite-dimensional simple real algebras}, Linear Algebra and its Applications \textbf{493} (2016), 164--182.
\bibitem{ShestakovZaicev} I. Shestakov, M. Zaicev, \textit{Polynomial identities of finite dimensional simple algebras} \textbf{2011}, 929--932.
\end{thebibliography}
\end{document}